# Space-time *hp*-finite elements for heat evolution in laser powder bed fusion additive manufacturing

Philipp Kopp[1] · Victor Calo[2] · Ernst Rank[1,3] · Stefan Kollmannsberger[1]

**Abstract**
The direct numerical simulation of metal additive manufacturing processes such as laser powder bed fusion is challenging due to the vast differences in spatial and temporal scales. Classical approaches based on locally refined finite elements combined with time-stepping schemes can only address the spatial multi-scale nature and provide only limited scaling potential for massively parallel computations. We address these shortcomings in a space-time Galerkin framework where the finite element interpolation also includes the temporal dimension. In this setting, we construct four-dimensional meshes that are locally refined towards the laser spot and allow for varying temporal accuracy depending on the position in space. By splitting the mesh into conforming time-slabs, we recover a stepwise solution to solve the space-time problem locally in time at this slab; additionally, we can choose time-slab sizes significantly larger than classical time-stepping schemes. As a result, we believe this setting to be well suited for large-scale parallelization. In our work, we use a continuous Galerkin–Petrov formulation of the nonlinear heat equation with an apparent heat capacity model to account for the phase change. We validate our approach by computing the AMB2018-02 benchmark, where we obtain an excellent agreement with the measured melt pool shape. Using the same setup, we demonstrate the performance potential of our approach by hatching a square area with a laser path length of about one meter.

## 1 Introduction

Laser powder bed fusion (LPBF) is an additive manufacturing technology that allows printing three-dimensional metal structures directly from a computer model. LPBF provides design flexibility and can be more efficient than other manufacturing techniques, for example, by reducing material waste and energy consumption. During the printing process, metal powder is added in a layerwise fashion and selectively melted by a high-power laser. The laser path is obtained by slicing the geometric model of the structure and hatching the interior areas. The path construction significantly impacts the microstructure and potential deviations from the geometric model (e.g., in zones with local overheating). Numerical simulations can identify these zones and significantly improve the quality in problematic areas by analyzing the manufacturing process and the structures without printing them. In particular, simulations can estimate quantities such as cooling rates that are not trivial to measure from experiments.

The major challenge in simulating LPBF processes is the multi-scale nature of the solution. While the laser spot and the melt pool sizes are well below 1 mm, the path length for manufacturing realistic structures can be in the range of kilometers. Similarly, the characteristic time scales range from micro-seconds around the laser spot to minutes or hours for the temperature evolution on the component level. This large span of scales renders detailed simulations of the melt pool dynamics unfeasible for simulating significant parts of the manufacturing process since relatively small fractions of the

✉ Philipp Kopp
philipp.kopp@tum.de

1. Chair of Computational Modeling and Simulation, Technische Universität München, Munich, Germany
2. Chair in Computational Geoscience & Applied Mathematics, School of Electrical Engineering, Computing & Mathematical Science, Curtin University, Perth, Australia
3. Institute for Advanced Study, Technische Universität München, Munich, Germany

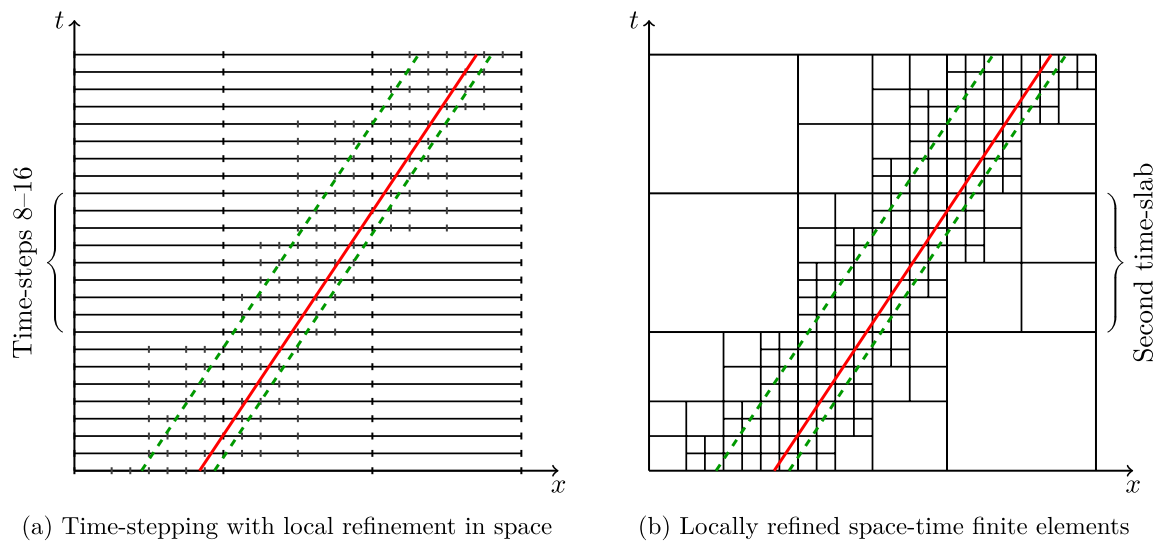

**Fig. 1** Discretization method comparison for a transient one-dimensional problem with a moving heat source (source path in red, refinement region in dashed green)

problem are often already expensive to compute. However, a reasonable approximation of the thermal history is often much cheaper to compute and can provide significant insight into the overall behavior. Based on these results, one can then select specific areas for a more detailed analysis which may also include the computation of the associated elastoplastic problem using classic techniques (using, e.g., [1]).

Our work focuses on the temperature evolution of laser-based additive manufacturing processes as it drives many phenomena and forms the basis for more detailed simulations. For example, the temperature evolution directly influences the microstructure formation (see, e.g., [2]) and can even be used to feed microstructural models [3]. Thermal models can predict melt pool shapes accurately (see, e.g., [4, 5] and references therein). In [6], the authors improve the quality consistency by comparing real-time measurements to a thermal simulation of the process. Other approaches, such as [7, 8], use data-driven models to predict the thermal history or mechanical properties based on the thermal history. Hence, developing efficient methods for simulating temperature evolution is an essential ongoing subject of research.

Nevertheless, even simple thermal models need to span many scales accurately and require particular strategies to render the computational requirements reasonable. Currently, a popular approach for direct simulations is to locally refine finite elements in space in combination with a time-stepping scheme. In particular, *hp*-finite elements were successfully applied in, e.g., [9, 10], where they reduce the computational effort by combining local refinement around the laser spot with high-order finite elements in the rest of the domain. The major shortcoming of this approach is that the time-step size is uniform in space and having sufficient accuracy around the laser spot leads to a needlessly accurate time integration on the rest of the domain. The resulting large number of time steps makes upscaling towards massively parallel simulations challenging as the communication overhead quickly dominates when computing individual time-steps in a distributed manner.

We address these challenges using locally refined space-time finite elements. We apply adaptive methods that capture the spatial and temporal multi-scale nature in a unified framework by considering time as a fourth dimension. Instead of solving the entire problem simultaneously, we split the space-time domain into several time-slabs that we compute sequentially. The size of each slab can be chosen depending on the target environment; for a single processor, one would choose much shorter time-slabs than for a supercomputer. Figure 1 schematically compares a classical time-stepping approach for solving a transient one-dimensional problem with 24 time-steps to a space-time discretization with local refinement towards a laser track that we compute across three time-slabs. Figure 2 shows a similar setup for two spatial dimensions where each time-slab is refined four times towards the laser source. While we believe that our way of simulating temperature evolution in LPBF is new, we want to mention [11–14] as related work on addressing the spatial and temporal multi-scale nature.

Many space-time finite element formulations for parabolic problems were developed over the last decades. Popular methods combine a continuous discretization in space with either a continuous Galerkin (cG) or discontinuous Galerkin (dG) approach in time. DG methods in time, such as discussed in [15, 16] allow for time marching schemes with non-matching discretizations as they enforce continuity on slab interfaces weakly. CG approaches tend to be more efficient for a given number of unknowns; therefore, they are suitable for local mesh refinement in space-time discretizations. Different formulations have been analyzed, such as [17, 18], and in particular [19], where they use truncated hierarchical B-Splines to develop an adaptive space-time formulation with higher continuity. Space-time finite element methods are also used to compute problems with moving boundaries [20, 21].

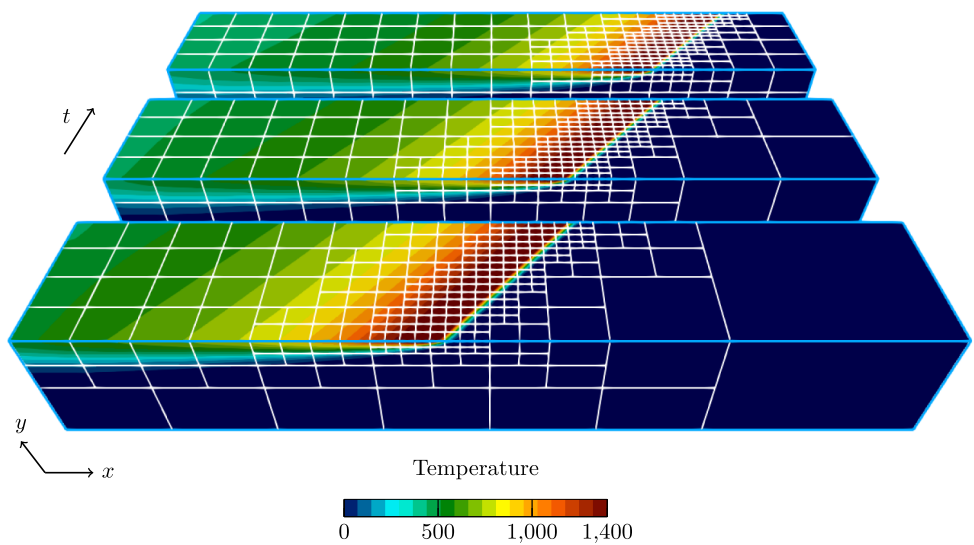

**Fig. 2** Three refined space-time slabs for a non-linear heat equation in two spatial dimensions

Our work picks up the formulation introduced in [22], where the test functions are the time derivative of functions from the trial space. Although the original formulation uses a tensor-product structure in space and time, the authors claimed that their method could form the basis for adaptive schemes; our present contribution confirms their claim. More recent publications, such as [23, 24], refer to our method as a continuous Galerkin–Petrov (cGP) formulation. The mixed derivative term of the form $\nabla \dot{w}$ resulting from testing with the time derivative is well defined as we use extruded spatial meshes with space tree refinement.

Developing continuous $hp$-finite element algorithms, especially in four dimensions, is a challenging task. Several approaches exist for three-dimensional problems, particularly the multi-level $hp$-method introduced in [25, 26]. In our previous work presented in [27], we developed a data-oriented approach for constructing multi-level $hp$-bases for higher dimensions by building on the core idea of using a hierarchical refinement where all levels can support basis functions. Our algorithms use basic adjacency information of cells in the refinement tree that require only a lightweight data structure. The simplicity of this version of the multi-level $hp$-method renders it an ideal choice for our space-time approach.

The paper is structured as follows. In Sect. 2, we discuss the nonlinear conductive model and its continuous Galerkin–Petrov formulation. We obtain a more accurate approximation of the melt pool geometry by including an apparent heat capacity term that accounts for the energy involved in the phase change. Sect. 3 recalls the construction of multi-level $hp$-bases and discusses the additional steps involved to ensure continuity between time-slabs. We then validate our approach in Sect. 4 by comparing the obtained melt-pool shape with the experimental results reported for the AMB2018-02 benchmark. We demonstrate the performance of our method by using this setup to compute a hatched square with a path length of one meter that computes in about seven hours on a single CPU. Finally, we conclude with some remarks and an outlook for future work in Sect. 5.

## 2 Formulation

We consider a heat equation with nonlinear coefficients on a spatial domain $\mathcal{S}$ and time interval $\mathcal{T} = [t_0, t_1]$:

$$\begin{aligned} c\dot{u} - \nabla \cdot (k\nabla u) &= f && \text{on } \Omega = \mathcal{S} \times \mathcal{T} \\ u &= u_0 && \text{on } \mathcal{S}, \text{ at } t = t_0 \\ u &= u_D && \text{on } \Gamma_D \\ n \cdot k\nabla u &= h && \text{on } \Gamma_N, \end{aligned} \quad (1)$$

where $u$ is the solution, $\dot{u}$ denotes its temporal derivative and $\nabla u$ its spatial gradient. We define the heat capacity $c = c(u)$, the heat conductivity $k = k(u)$, and $\Gamma_D \cup \Gamma_N = \partial \mathcal{S} \times \mathcal{T}$, such that $\Gamma_D \cap \Gamma_N = \emptyset$. We list the initial condition ($u_0$) and the Dirichlet boundary condition ($u_D$) separately, but we treat them similarly in our space-time formulation. Moreover, we expand the heat capacity into $c(u) = \rho c_s(u)$, where $\rho$ is the density and $c_s(u)$ is the specific heat capacity. In our examples, we use a volumetric heat source (Sect. 2.2) and homogeneous Neumann boundary conditions ($h = 0$). While not considered here, longer simulations of LPBF processes require including appropriate radiation and convection boundary conditions.

## 2.1 Weak form and linearization

To obtain the weak formulation of (1), we test with the time derivatives of the standard Bubnov–Galerkin test functions. The weak form of (1) then reads: Find $u \in u_b + W_0(\Omega)$, such that

$$\int_\Omega \dot{w}\, c\dot{u} + \nabla \dot{w} \cdot k\nabla u \, d\Omega = \int_\Omega \dot{w} f \, d\Omega + \int_{\Gamma_N} \dot{w}\, h \, d\Gamma_N \qquad \forall w \in W_0(\Omega), \quad (2)$$

where $u_b \in W(\Omega)$ satisfies the initial and Dirichlet conditions and $W_0(\Omega)$ restricts $W(\Omega)$ to functions that are zero on the boundary. The function space $W(\Omega)$ is the standard continuous finite element space $H^1(\Omega)$ with the additional condition that the gradient of the time derivative is defined. This condition on the time derivative excludes elements such as triangles in space-time where the Jacobian matrix is not diagonal (block diagonal in higher dimensions). Still, the continuity of the temporal derivative condition allows for locally refined rectangular elements with hanging nodes, such as shown in Fig. 1b. This type of Petrov-Galerkin formulation was introduced by [22]. Although they only consider tensor-products of spatial and temporal meshes, our results indicate that the formulation extends well to non-tensor-product spaces such as our *hp*-discretization.

We then choose a finite-dimensional subspace $W^h \subset W$, spanned by basis functions $N_i(x, t)$, to formulate the weak residual in terms of the discrete solution $u^h = N_i \hat{u}_i \in W^h$:

$$R_i(u^h) = \int_\Omega \dot{N}_i\, c\dot{u}^h + \nabla \dot{N}_i \cdot k\nabla u^h - \dot{N}_i f \, d\Omega - \int_{\Gamma_N} \dot{N}_i\, h \, d\Gamma_N, \quad (3)$$

with the linearization

$$T_{ij}(u^h) = \frac{\partial R_i}{\partial \hat{u}_j} = \int_\Omega \dot{N}_i \left(c\dot{N}_j + c'\dot{u}^h N_j\right) + \nabla \dot{N}_i \cdot \left(k\nabla N_j + k'\nabla u^h N_j\right) d\Omega, \quad (4)$$

where $c'$ and $k'$ are the derivatives of $c$ and $k$ with respect to temperature. With this, we obtain a nonlinear iteration scheme $\hat{u}^{k+1} = \hat{u}^k + \Delta \hat{u}^{k+1}$, where $\Delta \hat{u}^{k+1}$ is the solution of the linear equation system $T_{ij}^k \Delta \hat{u}_j^{k+1} = -R_i^k$. The superscripts in $R_i^k$ and $T_{ij}^k$ indicate that we evaluate (3) and (4) using $u_k^h = N_i \hat{u}_i^k$. We start with $\hat{u}^0 = 0$ and iterate until we obtain a reasonable reduction in the residual $\|R^k\| < \epsilon \|R^0\|$. We choose $\epsilon = 10^{-4}$ in our computations as smaller values did not improve the solution noticeably.

A consequence of (2) is that we obtain a discontinuous test space across time-slices where the finite element interpolation is $C^0$ continuous. Introducing such slices when meshing the space-time domain allows us to split the problem into time-slabs that we compute consecutively. Moreover, we eliminate a temporal test function on each slab when imposing the initial (Dirichlet) condition to recover linearly independent equations. See Sect. 3.3 for further discussion.

We use a standard Gauss-Legendre quadrature rule with $p + 1$ points to integrate the finite element systems in space. As the test functions have a maximum polynomial degree of $p - 1$ in time, it is even sufficient to use only $p$ quadrature points in time to integrate the products of shape functions accurately. However, in the vicinity of a concentrated heat source, we might need more quadrature points as the rapid temperature changes and the material's nonlinearities can make the integrals quite rough. This quadrature inadequacy is especially problematic when using an apparent heat capacity model to account for the phase change, as Sect. 2.3 discusses.

## 2.2 Heat source

Throughout this paper, we assume a Gaussian heat source of the form

$$g(x, y) = \frac{P\nu}{2\pi\sigma^2} \exp\left(-\frac{x^2 + y^2}{2\sigma^2}\right), \quad (5)$$

where $\sigma = D4\sigma/4$, $P$ is the laser power, and $\nu$ is the material's absorptivity. Instead of directly imposing $g$ as a heat flux boundary condition, we use a volumetric extension into the *z*-direction by an intensity function

$$I(z) = \frac{2}{\sqrt{2\pi}\sigma_z} \exp\left(-\frac{z^2}{2\sigma_z^2}\right), \quad (6)$$

for $z \leq 0$. Now, we define the source function *f* in terms of a given path $p(t)$ that determines the current position of the laser:

$$f(x, y, z, t) = g(x - p_x(t), y - p_y(t))\, I(z - p_z(t)). \quad (7)$$

Choosing a volumetric source has two advantages over imposing the heat flux on a surface. First, we account for the recoil-pressure-induced depressions around the laser spot that can form due to material evaporation. However, once these effects become dominant, we can no longer expect realistic results from a purely thermal model. Second, when modeling the metal powder as a continuum (although not done here), we can account for light rays naturally reaching some depth within the powder. We refer to [28] for an investigation on volumetric source shapes for temperature simulations of LPBF processes.

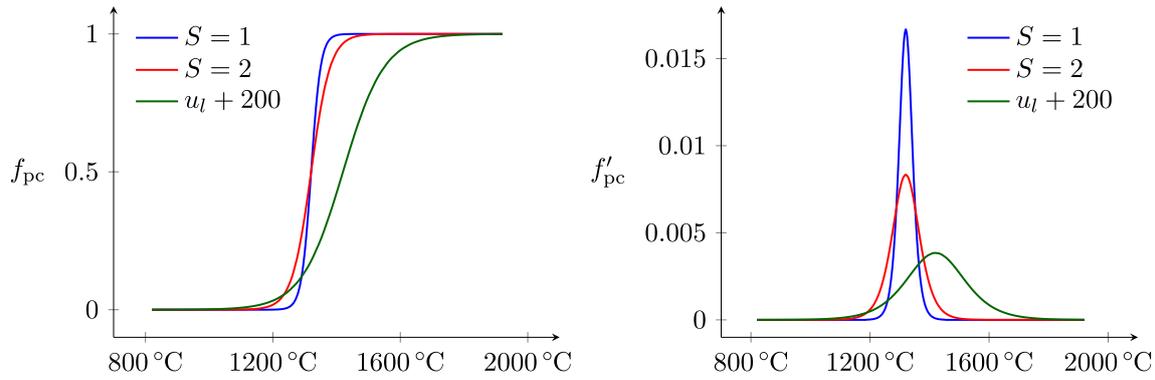

**Fig. 3** Phase change regularization between solid and liquid states

## 2.3 Phase change model

Without accounting for the energy involved in the phase change, we significantly overestimate the temperatures in the melt pool, which causes an increased heat energy dissipation around the laser spot, resulting in an underestimated melt pool length, and an overestimated melt pool width and depth. We improve this behavior by using an apparent heat capacity model that adds a latent heat contribution in a thin region around the phase change to the heat capacity. To this end, we define the following phase change function

$$f_{\text{pc}}(u) = \frac{1}{2}\left(\tanh\left(\frac{u - u_m}{u_\sigma}\right) + 1\right),$$

with

$$u_m = \frac{u_l + u_s}{2} \qquad u_\sigma = S\frac{u_l - u_s}{2}$$

that describes a smooth transition between the solid ($u_s$) and liquid ($u_l$) temperatures with a smoothness parameter $S$ (see Fig. 3). Using

$$f'_{\text{pc}}(u) = \frac{1}{2u_\sigma}\left(1 - \tanh\left(\frac{u - u_m}{u_\sigma}\right)^2\right)$$

$$f''_{\text{pc}}(u) = \frac{1}{u_\sigma^2}\tanh\left(\frac{u - u_m}{u_\sigma}\right)\left(\tanh\left(\frac{u - u_m}{u_\sigma}\right)^2 - 1\right),$$

we define the apparent heat capacity as

$$c(u) = \rho c_s(u) + \rho L f'_{\text{pc}}(u) \tag{8}$$

with its derivative

$$c'(u) = \rho c'_s(u) + \rho L f''_{\text{pc}}(u). \tag{9}$$

Here, $L$ represents the latent heat associated with the phase change from solid to liquid phases of the material. Our experience shows that the smoothness parameter $S$ mainly influences the length of the melt pool and the cooling rates below the melting temperature, which we exploit to fit the simulated melt pool length to the experimental data. We minimize the effect on the cooling rates by fixing $u_s$ while increasing $u_l$, which results in a shift of $u_m$ towards higher temperatures. Thus, the regularization influence is shifted to a range of temperatures that has little physical relevance due to the simplified model. These latent heat contributions are often neglected, as the nonlinearities render the solution process quite challenging. The quick and localized phase change induces a thin spike in the heat capacity $c(u)$ on the melt pool boundary that the discretization and the quadrature rule must capture. Besides increasing the number of quadrature points, we disregard the phase change in the initial iteration and then use a smooth regularization (e.g., $S = 10$) for the subsequent iteration. If necessary, we compute an improved update as $\hat{u}^{k+1} = \hat{u}^k + \beta \Delta \hat{u}^{k+1}$, where $\beta$ minimizes the residual norm and use a bisection type algorithm to prevent the iterations from diverging.

## 3 Discretization

Simulating real-world examples with $hp$-refinement in four dimensions requires an efficient basis function construction with manageable complexity. The multi-level $hp$-method introduced in [25, 26, 29] is an excellent framework for this task. In [27], we extend the framework to higher dimensions allowing us to build a space-time approach. We now summarize the main ideas of constructing multi-level $hp$-bases; see the references above for further details.

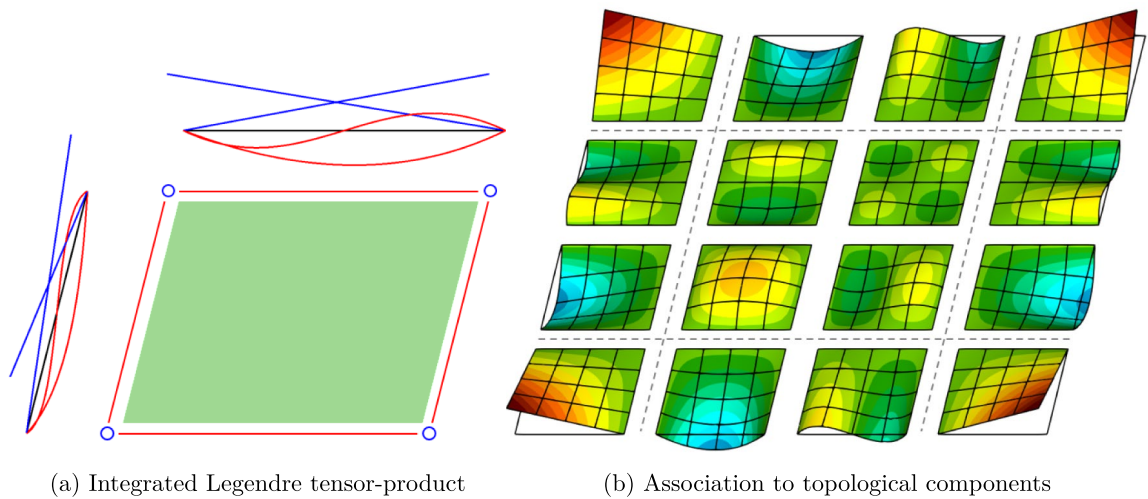

(a) Integrated Legendre tensor-product  (b) Association to topological components

**Fig. 4** Shape functions in the *p*-version of the finite element method

### 3.1 Multi-level *hp*-method

In the *p*-finite element method, tensor-products of integrated Legendre polynomials form higher-order $C^0$ continuous finite elements. The shape function hierarchy allows for varying polynomial degrees between elements as the functions on the side with a higher polynomial degree contain the functions on the other side. This hierarchical structure is commonly implemented by associating shape functions to the topological components of the element, resulting in a classification into nodal, edge, face, and volume modes. Figure 4 shows this association for a two-dimensional finite element. Additional advantages of using integrated Legendre polynomials are improved condition numbers (for Laplace type problems) and the possibility of defining a *trunk space* that removes certain functions from the *tensor-product space* without reducing the convergence order.

The multi-level *hp*-method allows to automatically refine a given base mesh by recursively building a space tree on top of it (quadtree in two dimensions, octree in three dimensions). Instead of replacing elements when refining, we keep the complete hierarchy and allow all cells (leaf or non-leaf) to support basis functions. This eliminates the problem of constraining hanging nodes at the cost of introducing slightly larger supports for basis functions in transition zones. Figure 5 compares a refinement by replacement to the hierarchical multi-level *hp* approach. As shown, shape functions on elements with finer neighbors are connected to the corresponding parent cell on the other side. To maintain global $C^0$ continuity, all overlay shape functions must be zero on internal boundaries.

We construct a multi-level *hp*-basis by first activating all basis functions associated with leaf cells' topological components and then deactivating the functions active on inner boundaries. Figure 6 shows the result of these two steps on a two-dimensional example. In both stages, we include all sub-components; when deactivating an edge, for example, we also deactivate the two nodes. Our extension to four dimensions expresses the same idea in terms of operations over element interfaces on the *tensor-product masks* associated with each cell. A tensor-product mask keeps track of active and inactive shape functions in the tensor-product of integrated Legendre polynomials. Thus, our data structure only needs to provide adjacency relations between cells, extending well to four dimensions. We skip the details of this formulation here and refer to the explanations given in [27].

In practice, on a hierarchically refined mesh, our implementation builds for every leaf cell a location map, an interface for evaluating the basis functions, and the indices of active face functions. Then, we use standard finite element technology that generally does not depend on the logic behind constructing the basis functions.

### 3.2 Refinement strategy

It is crucial to find strategies for automatic mesh refinement and polynomial degree selection to use the flexibility offered by an *hp*-framework to its full potential. General approaches use error estimators combined with some smoothness estimators to choose whether to refine in *h* or *p*. While these require little problem setup information, they are often expensive and complex to implement. However, in applications such as LPBF, we have an a priori knowledge about the laser path to tailor our mesh refinement strategies. We consider recent laser positions (e.g., 100 ms) to construct a function conceptually similar to a mesh density function that defines the target refinement depth. We start with a Gauss bell shape at the current position of the laser, where we have

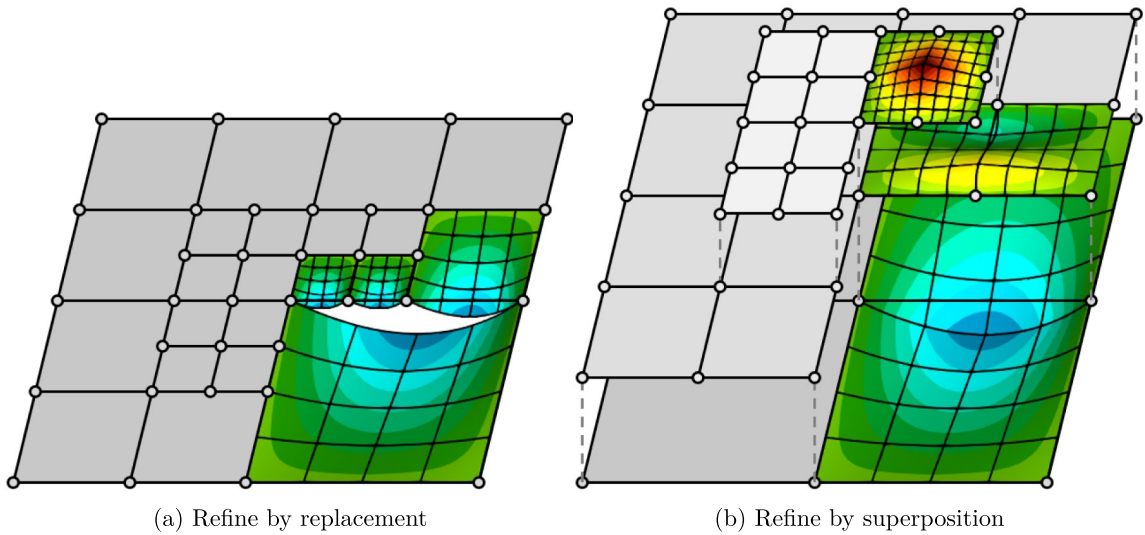

(a) Refine by replacement  (b) Refine by superposition

**Fig. 5** Comparison of *hp*-refinement strategies

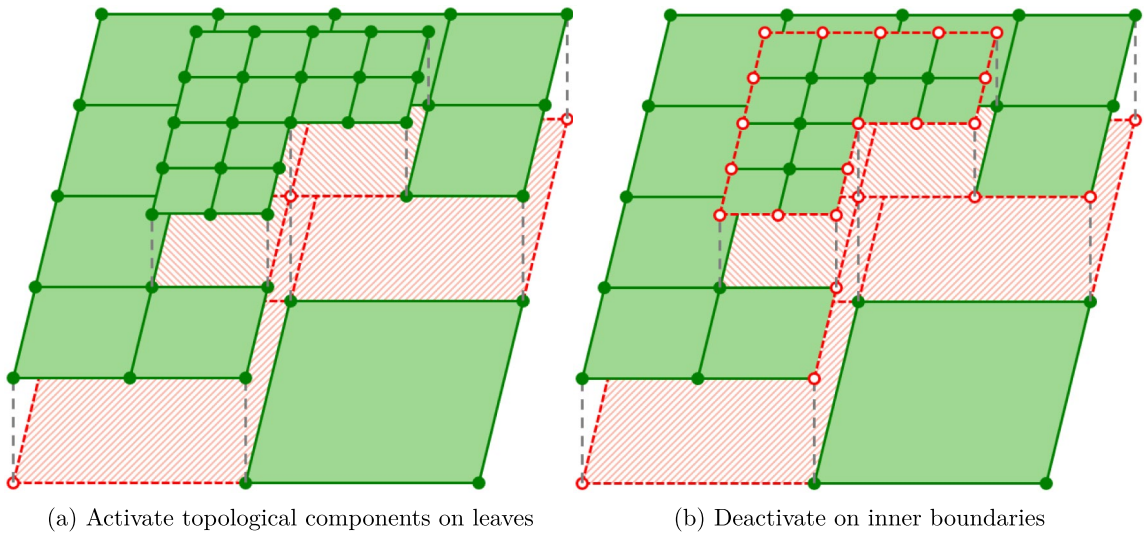

(a) Activate topological components on leaves  (b) Deactivate on inner boundaries

**Fig. 6** Multi-level *hp* construction rules. Deactivated components do not contribute their basis functions but may still support functions from lower dimensional components

a maximum refinement depth at the center that decreases as we move further away from the laser. Then, we transport this function in space along the previous laser path while decreasing the maximum refinement level and increasing the refinement width. This function transport gives us a spatial function for a given point in time with values indicating the target refinement level. We determine whether to refine an element by evaluating this function on a grid of points within the element (e.g., 5 to 7 per direction) and compare the maximum target refinement depth to the refinement level of the element. Figure 7 shows the maximum refinement depth $d_\tau$ and width $\sigma_\tau$ defined over the laser history that we later use in our examples.

We evaluate our refinement depth function $d(x, t)$ by considering the laser path in the time interval $[t - \tau_{\max}, t]$; the regions before and after do not influence the refinement. We first compute the closest point $x_i^p$ on each laser path segment $i$ and determine the time $t_i^p$ at which the laser was at $x_i^p$. Then, using the time delay $\Delta t_i = t - t_i^p$, we extract the maximum refinement depth and width from the functions $d_\tau(\Delta t)$ and $\sigma_\tau(\Delta t)$, respectively, as specified by, e.g., Figure 7. Using these, we compute $d(x, t)$ by evaluating a Gaussian function for each segment and taking the maximum value:

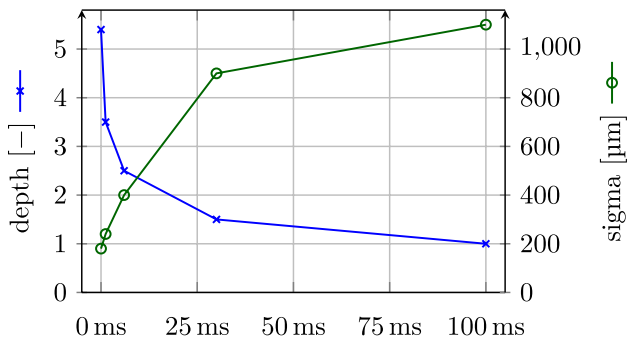

| delay [ms] | depth [−] | sigma [μm] |
|---|---|---|
| 0.0 | 5.4 | 180 |
| 1.2 | 3.5 | 240 |
| 6.0 | 2.5 | 400 |
| 30.0 | 1.5 | 900 |
| 100.0 | 1.0 | 1100 |

**Fig. 7** Parameters for refinement based on laser history

$$d(x,t) = \max_i \left( d_\tau \left(t - t_i^p\right) \exp\left( -\frac{\|x - x_i^p\|^2}{2\sigma_\tau^2(t - t_i^p)} \right) \right),$$

rounding the result to the closest integer. Figure 8 shows $d(x, t)$ for a square path with 6 mm width.

As our refinement is isotropic, the same refinement depth applies to all directions, including time. For simplicity, we choose the same polynomial degree for elements with equal refinement levels, which is not a restriction of the method. Finding better selection criteria for the polynomial distribution of individual elements will improve the simulation performance; future work will focus on this.

### 3.3 Separation into time-slabs

As Sect. 2 discusses, we test with the time derivatives of the trial functions. Let us momentarily consider a one-dimensional interpolation in time consisting of $C^0$ continuous quadratic finite elements, as Fig. 9a shows. The time derivatives (Fig. 9b) are linearly dependent as $\dot{N}_0 + \dot{N}_2 + \dot{N}_4 = 0$ and thus are not a valid test basis. As we eliminate the trial function $N_0$ by imposing the initial (Dirichlet) condition; thus, we eliminate the corresponding test function $\dot{N}_0$, which restores the linear independence to the system. For this, we use standard finite element routines for imposing Dirichlet boundary conditions during the assembly. Now, we can set up the monolithic system

$$\begin{pmatrix} \color{red}{K_{00}} & \color{red}{K_{01}} & \color{red}{K_{02}} & \color{red}{0} & \color{red}{0} \\ K_{10} & K_{11} & K_{12} & 0 & 0 \\ K_{20} & K_{21} & K_{22} & \color{red}{K_{23}} & \color{red}{K_{24}} \\ 0 & 0 & K_{32} & K_{33} & K_{34} \\ 0 & 0 & K_{42} & K_{43} & K_{44} \end{pmatrix} \begin{pmatrix} \hat{u}_0 \\ \hat{u}_1 \\ \hat{u}_2 \\ \hat{u}_3 \\ \hat{u}_4 \end{pmatrix} = \begin{pmatrix} \hat{f}_0 \\ \hat{f}_1 \\ \hat{f}_2 \\ \hat{f}_3 \\ \hat{f}_4 \end{pmatrix}$$

which after imposing the initial condition becomes

$$\begin{pmatrix} K_{11} & K_{12} & 0 & 0 \\ K_{21} & K_{22} & K_{23} & K_{24} \\ 0 & K_{32} & K_{33} & K_{34} \\ 0 & K_{42} & K_{43} & K_{44} \end{pmatrix} \begin{pmatrix} \hat{u}_1 \\ \hat{u}_2 \\ \hat{u}_3 \\ \hat{u}_4 \end{pmatrix} = \begin{pmatrix} \hat{f}_1 \\ \hat{f}_2 \\ \hat{f}_3 \\ \hat{f}_4 \end{pmatrix} - \hat{u}_0 \begin{pmatrix} K_{10} \\ K_{20} \\ 0 \\ 0 \end{pmatrix}.$$

We can solve this linear system with all four remaining unknowns at once, but due to the entries $K_{23}$ and $K_{24}$, we cannot first solve for $\hat{u}_1 - \hat{u}_2$ and then for $\hat{u}_3 - \hat{u}_4$. We can eliminate these values by removing the negative constant part on the right side of $\dot{N}_2$ (in Fig. 9b) by replacing $\dot{N}_2$ with $\dot{N}_2 + \dot{N}_4$. This does not change the solution as the resulting test functions span the same space as the original ones. Now we can solve the equivalent slab-wise scheme

$$\begin{pmatrix} K_{11} & K_{12} \\ K_{21} & K_{22} \end{pmatrix} \begin{pmatrix} \hat{u}_1 \\ \hat{u}_2 \end{pmatrix} = \begin{pmatrix} \hat{f}_1 \\ \hat{f}_2 \end{pmatrix} - \hat{u}_0 \begin{pmatrix} K_{10} \\ K_{20} \end{pmatrix}$$

$$\begin{pmatrix} K_{33} & K_{34} \\ K_{43} & K_{44} \end{pmatrix} \begin{pmatrix} \hat{u}_3 \\ \hat{u}_4 \end{pmatrix} = \begin{pmatrix} \hat{f}_3 \\ \hat{f}_4 \end{pmatrix} - \hat{u}_2 \begin{pmatrix} K_{32} \\ K_{42} \end{pmatrix},$$

where each new slab uses the interface unknowns from the previous one as initial conditions. Again, we can use standard finite element assembly routines to remove the respective rows in the element matrices and vectors when imposing initial conditions on each slab. Although our demonstration is one-dimensional, the process holds for our locally refined space-time meshes in four dimensions. In this case, eliminating the test functions of overlay meshes during the imposition of the initial condition like described above does not automatically remove all the negative parts of the test functions later in time, as the finer overlay meshes are not split into individual time-slabs. For example, on the first refinement level we obtain a basis in time similar to Fig. 9 for each root element, but only $N_0$ and $\dot{N}_0$ are eliminated when imposing the initial condition. Interestingly, we observe an improved conditioning during the sparse direct solution of the linear system when we set all negative constant parts of the test functions in time to zero during the assembly. This is similar to testing in time with Legendre polynomials of degree $p-1$ and not connecting them across element interfaces, instead of using the time derivatives of globally $C^0$ continuous integrated Legendre basis functions.

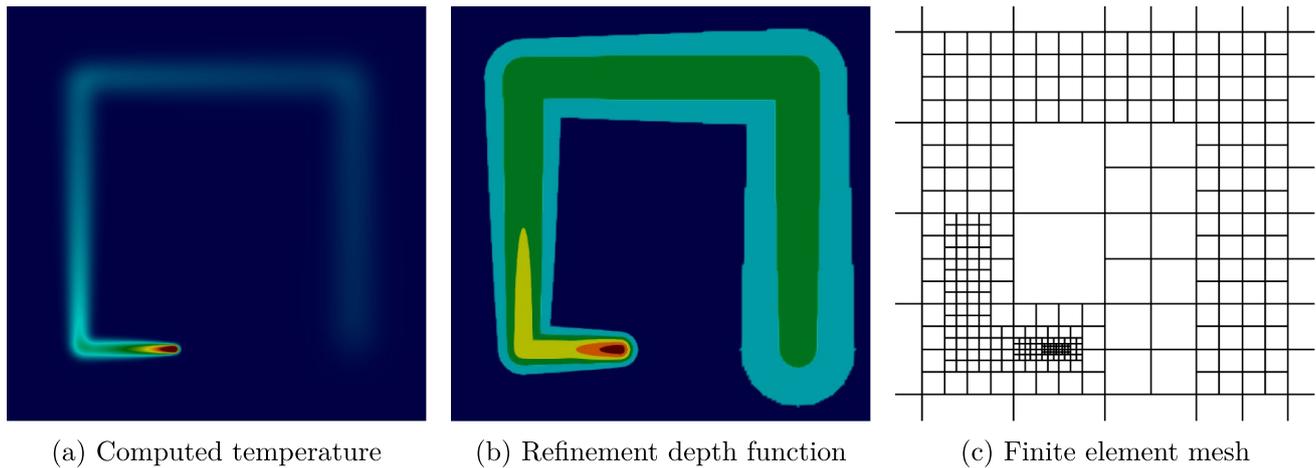

(a) Computed temperature  (b) Refinement depth function  (c) Finite element mesh

**Fig. 8** Mesh refinement using the parameters from Fig. 7 on a 6 mm by 6 mm square

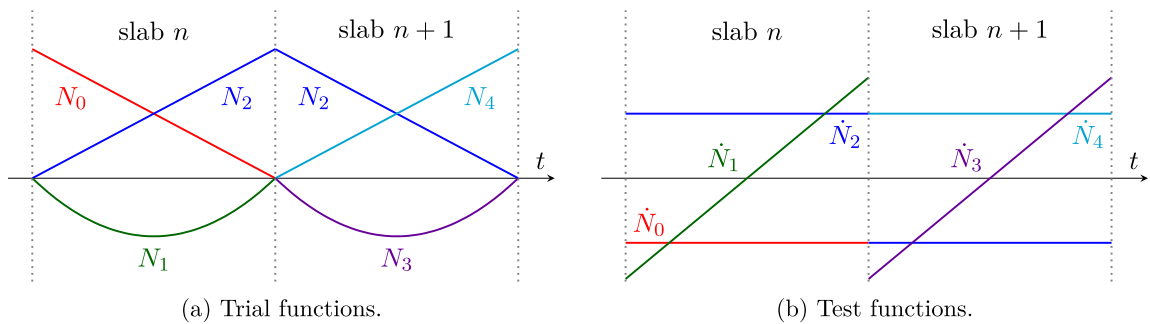

(a) Trial functions.  (b) Test functions.

**Fig. 9** Test and trial bases in time for two slabs

Another challenge in separating the solution into time-slabs is that hanging nodes are not constrained if we mesh each slab separately. Meshing the entire space-time domain at once is certainly not desirable. Still, we must at least know the refinement for the next slab and the previous one to construct compatible bases at slab interfaces. Figure 10 shows how we advance a one-dimensional solution in time by already meshing the second slab together with the first one. Then, when computing the second slab, we also mesh the third one, and only when computing the solution on the third slab do we discard the first mesh. This meshing procedure is possible as we define the refinement in advance independently of the solution. In other cases (e.g., when using an adaptive algorithm), we can either refine the next slab conformingly or use a dG formulation in time to enforce continuity across slab interfaces weakly.

## 4 Benchmark results

This section applies our methodology to the AMB2018-02 benchmark, which allows us to validate our approach and tune the model and discretization parameters. Using the same setup, we then compute a hatched square in Sect. 4.2. The source code is available at https://gitlab.com/hpfem/publications/2021_spacetime_am and its submodules under an open-source license. In addition to our space-time formulation, we provide time-stepping versions for all examples from this section to verify our implementation.

### 4.1 AMB2018-02

The AMB2018-02 benchmark [30] is commonly used to initially assess the performance of new simulation methods for LPBF processes. The setup features various laser configurations in single strokes on a block of IN625 alloy (24.08 mm × 24.82 mm × 3.18 mm). The length of each stroke is 14 mm with a duration of 17.5 ms, but the melt pool usually reaches a steady-state after at most 2 mm of travel distance. Experimental data from [31] shows the melt pool width and depth for each track as measured from the cross-section area within which the microstructure changed. The melt pool length and cooling rates are estimated using high refresh rate thermal imaging. Here, we focus on test case B on the AMMT setup and the measurements from [31, Table 3] for track number 3. Table 1 summarizes our simulation model parameters. The thermal properties of IN625 are

$$c_s(u) = \left(405 + \frac{247 \cdot u}{1000°C}\right) \left[\frac{J}{kg°C}\right] \qquad u \leq u_s$$

$$k(u) = \left(9.5 + \frac{15 \cdot u}{1000°C}\right) \left[\frac{W}{m°C}\right] \qquad u \leq u_s$$

for temperatures below the melting point and assumed constant for higher temperatures. Moreover, we use a constant density $\rho = 8440$ kg m$^{-3}$.

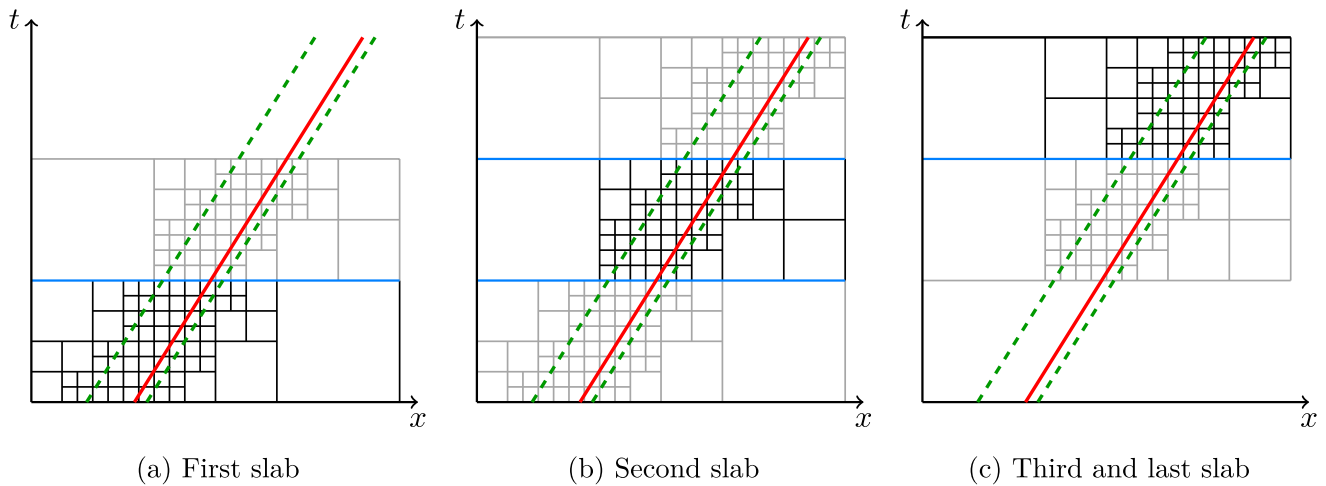

**Fig. 10** Ghost slabs in time with one space dimension

**Table 1** Model parameters (see Sects. 2.2 and 2.3 )

| Laser parameters | | Phase change parameters | |
|---|---|---|---|
| Speed ($v$) | 0.8 ms$^{-1}$ | Latent heat ($L$) | $2.8 \times 10^5$ J kg$^{-1}$ |
| Power ($P$) | 179.2 W | Solid temperature ($u_s$) | 1290°C |
| Absorptivity ($\nu$) | 0.32 | Liquid temperature ($u_l$) | 1350°C |
| Spot size ($D4\sigma$) | 170 μm | Initial temperature ($u_0$) | 25°C |
| Depth ($\sigma_z$) | $0.28 \cdot D4\sigma/4$ | Regularization ($S$) | 1, 2, 4 |

We calibrate our model in our first numerical experiments to reproduce the experimental data as close as possible. We estimate the melt pool dimensions from our simulations by extracting a contour surface at $u = u_s$ and measuring its bounding box's length, width, and depth. First, we identify suitable values for the absorptivity $\nu$ and the penetration depth $\sigma_z$ of the heat source as in (5) and (6) such that the width and depth match. Then, we set the phase change regularization for the melt pool length to match the experiments. We use this strategy since increasing the smoothness of the phase change mainly affects the melt pool length while it has a minor effect on its width and depth. We use a well-resolved discretization, see Table 2, to minimize the discretization influence. In our computations, we use linear polynomials in time combined with a trunk space for the spatial discretization (see, e.g., [27]). The additional $z$-factor in Table 2 multiplies the refinement width $\sigma_\tau$ in the $z$-direction to prevent unnecessary refinement in the region below the melt pool.

Table 3 shows the results for different phase change regularizations. We obtain very similar dimensions for $S = 4$ (equivalent to $S = 1$, $u_s = 1200$°C, and $u_l = 1440$°C) and with only increasing $u_l$ to 1550°C. However, Fig. 11 shows significant differences in temperatures below 1290°C if the phase change model is regularized towards lower temperatures. Figure 12 shows the time derivative of the temperature field that we obtain directly from the space-time finite element discretization. The cooling rates of around $2 \times 10^6$ °C s$^{-1}$ deviate significantly from the measured $1.08 \times 10^6$ °C s$^{-1}$, given by [31, Table 3]. However, in [31], the authors advise against using the cooling rate measurements due to motion blur and a limited calibration range. Therefore, the validation of cooling rates for our model remains an open task. Next, we coarsen the discretization as much as possible while still obtaining good estimates to showcase the benefits of our approach. Figure 13 shows the melt pool geometry and the temperature in the vicinity for the coarse discretization, see Table 2. With a simulation time of 17.5 ms and 8 time-slabs we obtain a duration of about 2.2 ms for one slab. The melt pool dimensions (last row of Table 3) are almost identical to ones of the well-resolved discretization. For simplicity, we over-integrate all elements with $p + 2$ quadrature points in space and $p + 1$ quadrature point in time to capture the phase change accurately. We solve the linear systems using Intel's Pardiso sparse direct solver.

**Table 2** Fine and coarse discretization parameters

| Base mesh | Fine discretization | | | | Coarse discretization | | | |
|---|---|---|---|---|---|---|---|---|
| | 64 × 64 × 9 × 32 slabs | | | | 12 × 12 × 3 × 8 slabs | | | |
| Refinement | $\tau$ [ms] | $d_\tau$ | $\sigma_\tau$ [μm] | z-factor | $\tau$ [ms] | $d_\tau$ | $\sigma_\tau$ [μm] | z-factor |
| | 0 | 6 | 100 | 0.5 | 0 | 5.4 | 180 | 0.5 |
| | 0.08 | 5.3 | 120 | 0.5 | 1.2 | 3.5 | 240 | 0.5 |
| | 0.47 | 4.5 | 150 | 0.5 | 6 | 2.5 | 400 | 0.8 |
| | 1.2 | 3.5 | 160 | 0.8 | 30 | 1.5 | 900 | 1 |
| | 4 | 2.5 | 200 | 1 | 100 | 1 | 1100 | 1 |
| | 16 | 0.5 | 240 | 1 | | | | |
| Polynomial degrees in space | 1 for level 0 and 2 for levels 1 to 6 | | | | 2, 2, 4, 4, 3 and 3 for levels 0 to 5 | | | |

**Table 3** Melt pool dimensions for different model parameters in comparison to experimental data

| Result | Length [μm] | Width [μm] | Depth [μm] |
|---|---|---|---|
| Measurements | 359 ($\sigma = 20$) | 132 ($\sigma = 2$) | 36 ($\sigma = 0.9$) |
| No latent heat | 301 | 138 | 39.4 |
| $S = 1$ | 396 | 129 | 34.8 |
| $S = 2$ | 381 | 129 | 34.8 |
| $S = 4$ | 356 | 129 | 34.8 |
| $S = 8$ | 328 | 130 | 34.8 |
| $u_l + 200$ | 354 | 131 | 35.4 |
| $u_l + 200$ (coarse) | 353 | 132 | 35.3 |

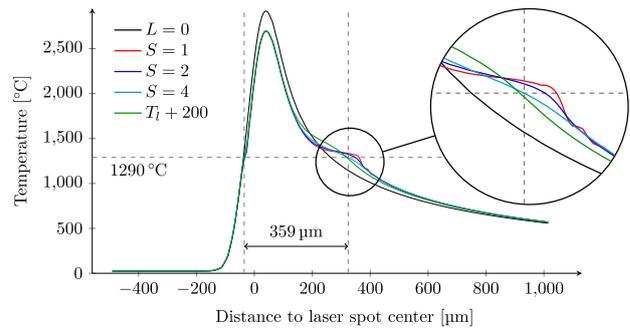

**Fig. 11** Temperature in laser travel direction for different phase change parameters

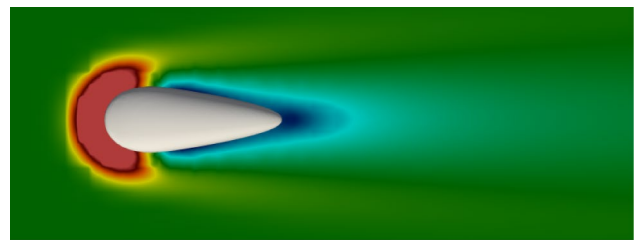

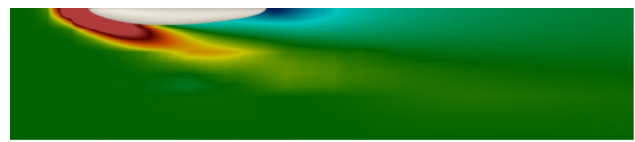

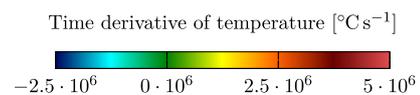

**Fig. 12** Time derivative of the temperature field using our fine discretization

### 4.2 Hatched square

In this section, we use the same setup of the previous section for the AMB2018-02 benchmark and hatch a square area with a 10 mm side length. As Fig. 14 shows, the laser path first follows the boundary and then fills the interior with a hatch distance of 100 μm, resulting in a path length of about 102 cm. The entire process takes about 1.28 s, which we extend to a total simulation time of 3 s. During this cooldown period, our discretization is automatically coarsened in space and in time as a result of the formulation and the way we construct our meshes, as Sect. 3.2 discusses. Each time-slab contains one base element in time with a duration of 2.4 ms, which is slightly longer than in the previous example. The number of unknowns per slab initially averages around 50 thousand and drops to around 2400 for slabs in the cooldown period.

Figure 15 shows the solution and the melt pool dimensions for two time-slices; on the left side, the laser approaches the top left corner, and on the right side, it has just reached it. We observe significant heating of the plate leading to more than a 50% predicted increase in melt pool

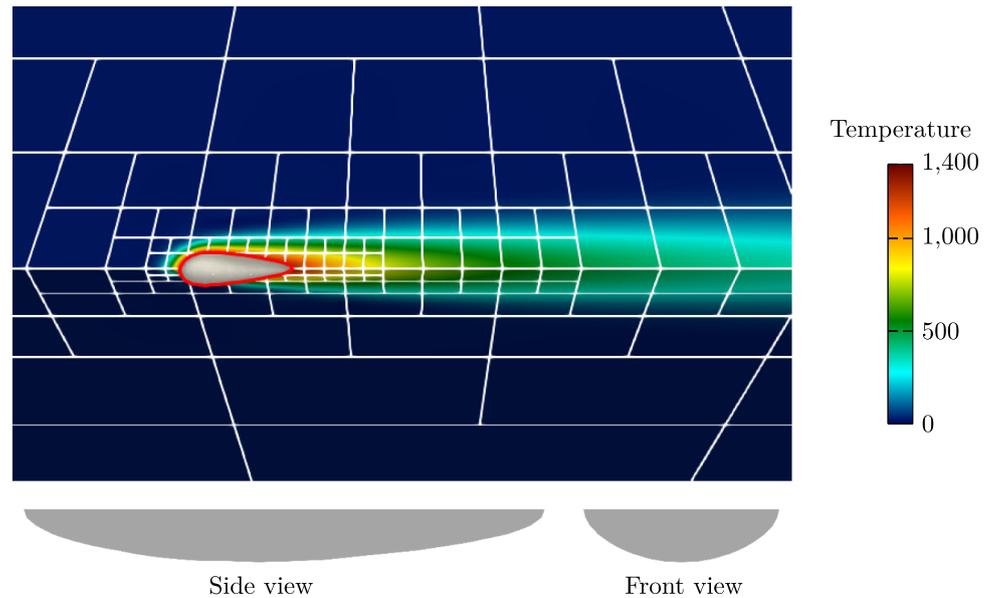

**Fig. 13** Solution and contour surface at $u_s = 1290°C$ for coarse discretization with $u_l = 1550°C$

width and depth for the second case. Hence we identify this spot as a potential source for defects.

The simulation runs on a single Intel Xeon Gold 6230 CPU with 20 cores running at 2.1 GHz in about 7 h. Compared to the single-threaded execution, our parallel version is about 8 times faster, a good result considering that the CPU's turbo boost frequency is 3.9 GHz, which from experience results in a maximum speedup of 11 to 12 times. While we reach this speedup in the assembly of the linear systems, the Pardiso sparse direct solver does not scale optimally in our examples. Moreover, we over-integrate again with $p + 2$ points in space and $p + 1$ points in time. We can improve this by identifying the elements around the laser source and increasing the number of quadrature points there only.

## 5 Conclusion

We present a space-time finite element approach for simulating heat evolution in metal additive manufacturing fusion processes. We can control the size of the resulting systems by choosing the desired number of elements in time, which increases the number of unknowns per slab compared with conventional time-stepping schemes. This control is essential for the scalability of parallel implementations as very small time-steps prevent good scaling due to the communication overhead. We obtain good speedups using shared memory parallelism on a single CPU. We expect this to transfer to large-scale distributed memory parallelism if we can find suitable linear solvers. Therefore, potential directions for future research are BDDC preconditioners, such as the one used in [32], or symmetric formulations, such as presented in [33].

Combining space-time finite elements with the multi-level hp-method allows us to introduce local mesh refinement in four dimensions to capture the multi-scale nature in both space and time. As a result, our approach allows for a much lower accuracy away from the heat source by either reducing the element duration (element size) or decreasing the polynomial degree in time. Although we did not yet observe a significant speedup over comparable time-stepping approaches with hp-refinement only in space, we expect this to change when computing realistic problems where the difference in scales is more significant. One reason for this behavior is a substantial overlap between the supports of basis functions in the refined region; therefore, finding strategies for reducing the size of the supports or even lowering the continuity may improve the performance even further.

We also show how to use the laser path to construct a tailored refinement strategy that manages to capture the dynamics of the process very well while not depending on the solution. This approach can form the basis for more advanced adaptive strategies based on error estimators to automate the refinement process further. Especially a more intelligent selection of the polynomial degrees would be desirable (as opposed to choosing one value per refinement level). A challenge is assessing the cost for h- and p-refinement, which is not always intuitive as it depends on the mesh topology. For example, increasing the polynomial degree in a coarse transition element makes the smaller elements in the vicinity significantly more expensive.

Finally, we demonstrate how a simple single field thermal model with a latent heat contribution can accurately reproduce the experimental data on the melt pool dimensions in the context of the AMB2018-02 benchmark. In the future, we also plan to validate geometrically more complicated laser paths and extend the model to account for the reduced conductivity in the metal powder that has not yet been melted. To simulate more extended periods, we also must include the effects of radiative and convective cooling on the surface, which we did not consider here. With this, we can start optimizing manufacturing processes where we may benefit from extracting cooling rates directly from our space-time finite element formulation.

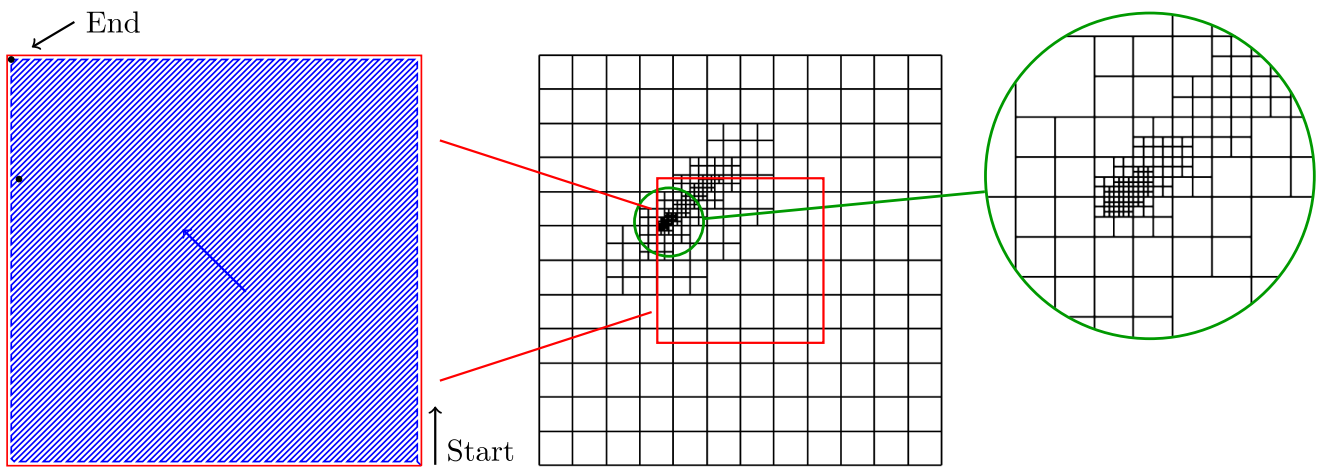

**Fig. 14** Laser path (left) and example discretization (center and right) at $t = 1.2168$ s for hatched square with 10 mm side length and 100 µm hatch width

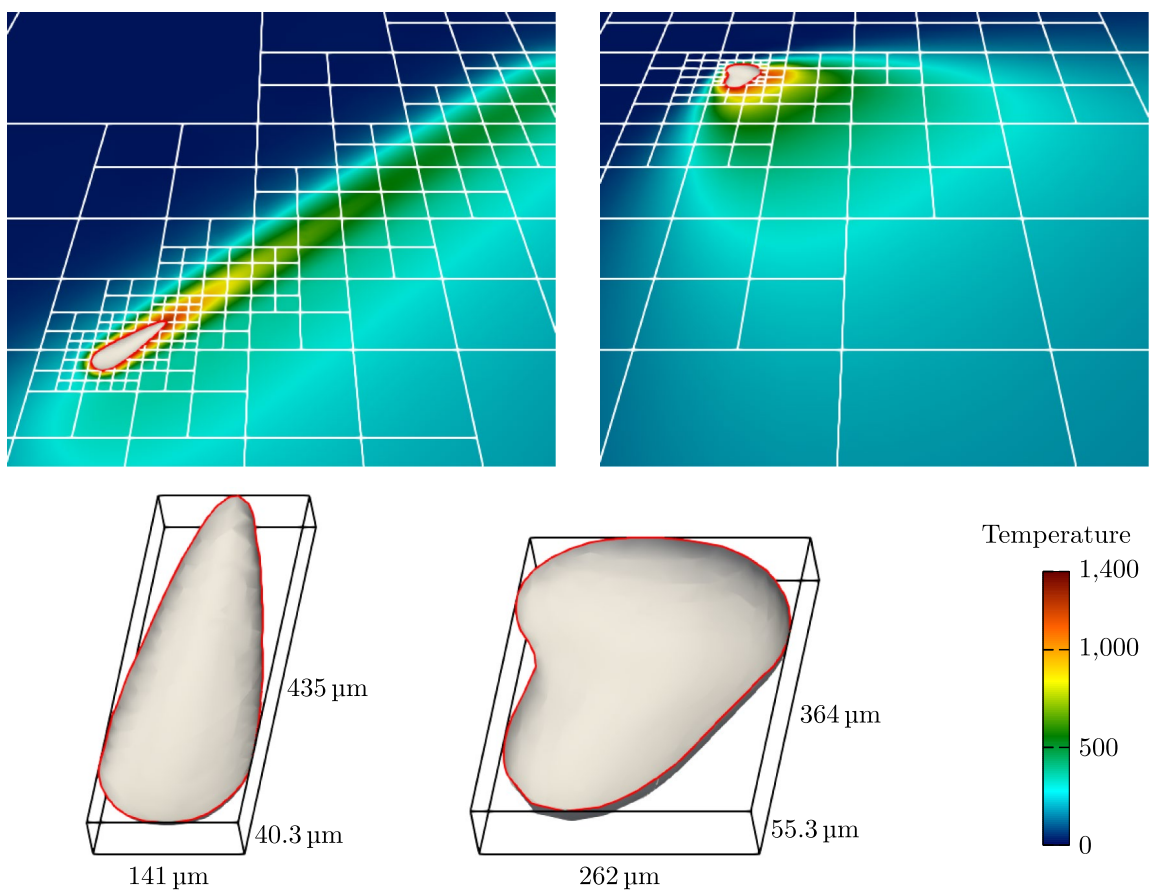

**Fig. 15** Temperature and melt pool geometry for $t = 1.2168$ s (left) and $t = 1.2786$ s (right)

**Acknowledgements** We gratefully acknowledge the support of Deutsche Forschungsgemeinschaft (DFG) through the grant KO 4570/2-1. This publication was also made possible in part by the CSIRO Professorial Chair in Computational Geoscience at Curtin University. This project has received funding from the European Union's Horizon 2020 research and innovation programme under the Marie Sklodowska-Curie grant agreement No 777778 (MATHROCKS). The Curtin Corrosion Centre kindly provides ongoing support.

## Declarations

**Conflict of interest** No potential conflict of interest was reported by the authors.

## References

1. Chiumenti M, Cervera M, Moreira Filho C, Barbat G (2020) Stress, strain and dissipation accurate 3-field formulation for inelastic isochoric deformation. Finite Elements in Analysis and Design, vol. 192, 12
2. Gu D, Shi Q, Lin K, Xi L (2018) Microstructure and performance evolution and underlying thermal mechanisms of Ni-based parts fabricated by selective laser melting. Additive Manuf 22:05
3. Nitzler J, Meier C, Müller K, Wall W, Hodge N (2021) A novel physics-based and data-supported microstructure model for part-scale simulation of laser powder bed fusion of Ti-6Al-4V. Adv Model Simul Eng Sci 8:12
4. Kollmannsberger S, Carraturo M, Reali A, Auricchio F (2019) Accurate prediction of melt pool shapes in laser powder bed fusion by the non-linear temperature equation including phase changes - isotropic versus anisotropic conductivity. Integrating Materials and Manufacturing Innovation, vol. 8, p. 167-177, 03
5. Paulson N, Gould B, Wolff S, Stan M, Greco A (2020) Correlations between thermal history and keyhole porosity in laser powder bed fusion. Additive Manufacturing, vol. 34, p. 101213, 04
6. Zhong Q, Tian X, Huo C (2021) Using feedback control of thermal history to improve quality consistency of parts fabricated via large-scale powder bed fusion. Additive Manufacturing, 03
7. Roy M, Wodo O (2020) Data-driven modeling of thermal history in additive manufacturing. Additive Manuf. vol. 32, p. 101017, 03
8. Xie X, Bennett J, Saha S, Lu Y, Cao J, Liu W, Gan Z (2021) Mechanistic data-driven prediction of as-built mechanical properties in metal additive manufacturing. npj Computational Materials. vol. 7, p. 86, 06
9. Kollmannsberger S, Özcan A, Carraturo M, Zander N, Rank E (2017) A hierarchical computational model for moving thermal loads and phase changes with applications to selective laser melting. Comput Math Appl 75:12
10. Kollmannsberger S, Kopp P (2021) On accurate time integration for temperature evolutions in additive manufacturing. GAMM-Mitteilungen 44:11
11. Soldner D, Mergheim J (2019) Thermal modelling of selective beam melting processes using heterogeneous time step sizes. Computers & Mathematics with Applications, vol. 78, no. 7, pp. 2183–2196. Simulation for Additive Manufacturing
12. Hodge N (2021) Towards improved speed and accuracy of laser powder bed fusion simulations via representation of multiple time scales. Additive Manuf 37:101600
13. Cheng L, Wagner GJ (2021) An optimally-coupled multi-time stepping method for transient heat conduction simulation for additive manufacturing. Comput Methods Appl Mech Eng 381:113825
14. Viguerie A, Carraturo M, Reali A, Auricchio F (2021) A spatiotemporal two-level method for high-fidelity thermal analysis of laserpowder bed fusion. 10. arXiv:2110.12932
15. Jamet P (1978) Galerkin-type approximations which are discontinuous in time for parabolic equations in a variable domain. Siam Journal on Numerical Analysis - SIAM J NUMER ANAL, vol. 15, pp. 912–928, 10
16. Eriksson K, Johnson C, Thomée V (1985) Time discretization of parabolic problems by the discontinuous galerkin method. J Multivariate Anal MA 19:01
17. Steinbach O (2015) Space-time finite element methods for parabolic problems. Comput Methods Appl Math 15:01
18. Devaud D, Schwab C (2018) Space-time $hp$-approximation of parabolic equations. Calcolo 55:09
19. Langer U, Matculevich S, Repin S (2019) 5. Adaptive space-time isogeometric analysis for parabolic evolution problems, pp. 141–184. De Gruyter, 09
20. Tezduyar T, Behr M, Liou J (1992) A new strategy for finite element computations involving moving boundaries and interfaces—the deforming-spatial-domain/space-time procedure: I. the concept and the preliminary numerical tests. Computer Methods in Applied Mechanics and Engineering, vol. 94, pp. 339–351, 02
21. Karyofylli V, Wendling L, Make M, Hosters N, Behr M (2019) Simplex space-time meshes in thermally coupled two-phase flow simulations of mold filling. Comput Fluids 192:104261
22. Aziz A, Monk P (1989) Continuous finite elements in space and time for the heat equation. Mathematics of Computation, vol. 52, pp. 255–274, 04
23. Schieweck F (2010) A -stable discontinuous galerkin–petrov time discretization of higher order. J Numer Math J NUMER MATH. vol. 18, pp. 25–57, 04
24. Hussain S, Schieweck F, Turek S (2011) Higher order galerkin time discretizations and fast multigrid solvers for the heat equation. J Numer Math 19:05
25. Zander N, Bog T, Kollmannsberger S, Schillinger D, Rank E (2015) Multi-level $hp$-adaptivity: high-order mesh adaptivity without the difficulties of constraining hanging nodes. Comput Mech 55:499–517
26. Zander N, Bog T, Elhaddad M, Frischmann F, Kollmannsberger S, Rank E (2016) The multi-level $hp$-method for three-dimensional problems: Dynamically changing high-order mesh refinement with arbitrary hanging nodes. Comput Methods Appl Mech Eng 310:252–277
27. Kopp P, Rank E, Calo V, Kollmannsberger S (2021) Efficient multi-level $hp$-finite elements in arbitrary dimensions. 06. arXiv:2106.08214
28. Zhang Z, Huang Y, Rani Kasinathan A, Imani Shahabad S, Ali U, Mahmoodkhani Y, Toyserkani E (2019) 3-dimensional heat transfer modeling for laser powder-bed fusion additive manufacturing with volumetric heat sources based on varied thermal conductivity and absorptivity. Optics & Laser Technology, vol. 109, pp. 297–312
29. Zander N, Bériot H, Hoff C, Kodl P, Demkowicz L (2022) Anisotropic multi-level hp-refinement for quadrilateral and triangular meshes. Finite Elements Anal Design 203:103700
30. AMB2018-02 Description. https://www.nist.gov/ambench/amb2018-02-description. Accessed: 2021-11-04
31. Lane B, Heigel J, Ricker R, Zhirnov I, Khromschenko V, Weaver J, Phan T, Stoudt M, Mekhontsev S, Levine L (2020) Measurements of melt pool geometry and cooling rates of individual laser traces on IN625 bare plates. Integrating Materials Manuf Innovation 9:02
32. Langer U, Yang H (2020) BDDC Preconditioners for a Space-time Finite Element Discretization of Parabolic Problems. pp. 367–374. International Conference on Domain Decomposition Methods, 10
33. Führer T, Karkulik M (2021) Space–time least-squares finite elements for parabolic equations. Computers & Mathematics with Applications, vol. 92, pp. 27–36, 06